\documentclass[12pt]{article}
\usepackage{amssymb}
\usepackage{amssymb}
\usepackage{amsmath, amscd}
\usepackage{color}
\usepackage{latexsym}
\usepackage{graphicx}
\usepackage{changes}
\usepackage{url}
\newtheorem{theorem}{Theorem}[section]
\newtheorem{lemma}[theorem]{Lemma}
\newtheorem{corollary}[theorem]{Corollary}
\newtheorem{definition}[theorem]{Definition}

\newcommand{\qed}{\hfill $\Box$ }

\newcommand{\proof}{\noindent{\bf Proof.}\ \ }
\begin{document}

\title{\Large {\bf  Resonance graphs of plane bipartite graphs as daisy cubes}}

\maketitle

\begin{center}
{\large \bf Simon Brezovnik$^{a,b}$, Zhongyuan Che$^{c}$,
Niko Tratnik$^{b,d}$,\\ Petra \v Zigert Pleter\v sek$^{d,e}$
}
\end{center}
\bigskip\bigskip

\baselineskip=0.20in

\smallskip

$^a$ {\it University of Ljubljana, Faculty of Mechanical Engineering, Slovenia} \\

$^b$ {\it Institute of Mathematics, Physics and Mechanics, Ljubljana, Slovenia} \\

$^c$ {\it Department of Mathematics, Penn State University, Beaver Campus, Monaca, USA} \\

$^d$ {\it University of Maribor, Faculty of Natural Sciences and Mathematics, Slovenia} \\
\medskip

$^e$ {\it University of Maribor, Faculty of Chemistry and Chemical Engineering, Slovenia}\\
\medskip

\begin{center}
{\tt simon.brezovnik@fs.uni-lj.si, zxc10@psu.edu,\\ niko.tratnik@um.si, petra.zigert@um.si}
\end{center}

\begin{abstract}
We characterize plane bipartite graphs whose resonance graphs are daisy cubes, and therefore
generalize related results on resonance graphs of benzenoid graphs, catacondensed even ring systems, 
as well as 2-connected outerplane bipartite graphs.
Firstly, we prove that if $G$ is a plane elementary bipartite graph other than $K_2$, 
then the resonance graph of $G$ is a daisy cube if and only if the Fries number of $G$ equals the number of finite faces of $G$. 
Next, we extend the above characterization from plane elementary bipartite graphs to plane bipartite graphs and show that 
the resonance graph of a plane bipartite graph $G$ is a daisy cube if and only if $G$ is weakly elementary bipartite
such that each of its elementary component $G_i$ other than $K_2$ holds the property that the Fries number of $G_i$ 
equals the number of finite faces of $G_i$.
Along the way, we provide a structural characterization 
for a plane elementary bipartite graph whose resonance graph is a daisy cube,
and show that a Cartesian product graph is a daisy cube if and only if all of its nontrivial factors are daisy cubes.


\end{abstract}

\section{Introduction}

Resonance graphs illustrate the interconnections among perfect matchings, known as Kekulé structures in the realm of chemistry. 
These graphs were  first introduced independently by chemists El-Basil \cite{el-basil-93/1,el-basil-93/2} 
and Gründler \cite{grundler-82}. While the research on resonance graphs initially revolved around benzenoid graphs \cite{zhgu-88},
 the concept was subsequently extended to encompass plane bipartite graphs,
coined as $Z$-transformation graphs by mathematicians Zhang, Guo, and Chen \cite{zhgu-88},
resulting in a broader application of this framework. 
Some results on the role of hypercubes in  resonance graphs of benzenoid graphs can be found,
 for example, in \cite{SKG06,SKVZ09,TZ12}. 
 On the other hand, structural properties of  resonance graphs of plane bipartite graphs are gathered in the survey paper \cite{Z06}, 
 see also  \cite{LZ03,TV12,ZLS08,ZZ00,ZZY04}. 
 For some novel results on this topic see \cite{br-tr-zi-2, BCTZ23, C18, C19, C21}. 
 Recently, the concept of  resonance graphs was extended further to graphs embedded on closed surfaces \cite{TD23}.

Daisy cubes introduced in \cite{KM19} are a subfamily of partial cubes which contains Fibonacci cubes and Lucas cubes.
Some recent results on the characterizations and proper embeddings of daisy cubes can be found in \cite{T20,V21}. 
It was proved in \cite{KZ05} that Fibonacci cubes are resonance graphs of fibonaccenes, i.e., zigzag hexagonal chains.
In \cite{ZOY09}, Zhang et al. characterized plane bipartite graphs whose resonance graphs are Fibonacci cubes, 
and proved that Lucas cubes cannot be resonance graphs.
They also showed that Fibonacci cubes cannot be Cartesian products of nontrivial graphs.

A \textit{catacondensed benzenoid graph} is a 2-connected outerplane bipartite graph embedded in 
a hexagonal lattice such that every interior region is bounded by a unit hexagon. 
Let $h$ be a hexagon of a catacondensed benzenoid graph.
Then $h$ is an \textit{angularly connected hexagon}  if 
$h$ contains exactly two vertices of degree 2 and such two vertices are adjacent, and a \textit{linearly connected hexagon} 
if $h$ contains exactly two vertices of degree 2 and such two vertices are nonadjacent. 
By definitions, we can see that any angularly (respectively, linearly) connected hexagon
has an common edge with two other hexagons in a catacondensed benzenoid graph.
A catacondensed benzenoid graph without any linearly connected hexagons is called a \textit{kinky benzenoid graph}.
The connection between resonance graphs and daisy cubes was firstly investigated in \cite{Z18}, 
where it was shown that the resonance graph of a kinky benzenoid graph is a daisy cube. 
A \textit{catacondensed even ring system} is a 2-connected outerplane bipartite  graph whose each vertex has degree 
at most 3. Catacondensed even ring systems and 2-connected outerplane bipartite graphs
whose resonance graphs are daisy cubes were characterized in \cite{br-tr-zi-1} and \cite{br-tr-zi-2}  respectively,
by generalizing  angularly connected hexagons and  linearly connected hexagons
of a catacondensed benzenoid graph to equivalent concepts of a 2-connected outerplane bipartite graph.

The aim of this paper is to characterize plane bipartite graphs whose resonance graphs are daisy cubes. 
In Section \ref{sec2},
we give necessary concepts and known results. In Section \ref{sec3}, 
we focus on plane elementary bipartite graphs and
provide a list of equivalent characterizations for this family of graphs whose resonance graphs are daisy cubes.
In Section \ref{sec4}, 
we first show that a Cartesian product graph is a daisy cube if and only if each of its nontrivial factors is a daisy cube. 
We then apply this result to obtain characterizations of plane bipartite graphs  whose resonance graphs are daisy cubes, 
by considering their elementary components.

\section{Preliminaries} \label{sec2}

For a graph $G$,  $V(G)$ is the set of vertices, and  $E(G)$ is the set of edges of $G$. 
 A subgraph of $G$ induced by a vertex subset $X \subseteq V(G)$ is denoted as $\langle X \rangle$. 
The degree of a vertex $u \in V(G)$ is denoted by $\deg_G(u)$. 
For vertices $u,v \in V(G)$, let $d_G(u,v)$ denote the distance between $u$ and $v$,
which is the length of a shortest path  between two vertices;
and let  $I_G(u,v)$ denote the set of  all vertices  that are on shortest paths between $u$ and $v$ in $G$.
Assume that $H$ is a connected  induced subgraph of $G$. Then 
$H$ is called an  \textit{isometric subgraph} of $G$ if  $d_H(u,v)=d_G(u,v)$ for any two vertices $u,v\in V(H)$.
A connected graph $G$ is called a \textit{median graph} 
if  {$|I_G(u,v) \cap I_G(u,w) \cap I_G(v, w)|=1$} for every triple of vertices $u, v, w$ of $G$.

\subsection{Daisy Cubes}\label{S2-DaisyCubes}

A \textit{Cartesian product} of graphs $G_1, G_2, \ldots, G_t$, 
denoted by $\Box_{i=1}^t G_i=G_1 \Box G_2 \Box \cdots \Box G_t$, where $t \ge 2$,
is a graph with the vertex set $V(G_1) \times V(G_2) \times \cdots \times V(G_t)$, 
and two vertices $a=(a_1,a_2,\ldots, a_t)$, $b = (b_1, b_2, \ldots, b_t)$ 
are adjacent if  there exists exactly one $j \in \{ 1, \ldots, t \}$ such that $a_jb_j \in E(G_j)$, 
and $a_i=b_i$ for all $i \in \{ 1, \ldots, t \} \setminus \{ j \}$. 
It is well known that a Cartesian product $\Box_{i=1}^t G_i$ is connected 
if and only if each factor $G_i$ is connected  for all $1 \le i \le t$. 
Moreover, for any two vertices $a,b$ of $\Box_{i=1}^t G_i$, the following distance formula holds true \cite{HIK11}:
\begin{equation} \label{en_dis}
d_{\Box_{i=1}^t G_i}(a,b) = \sum_{i=1}^t d_{G_i}(a_i,b_i).
\end{equation}

A hypercube $Q_0$ of dimension 0 is the one-vertex graph $K_1$, 
a hypercube $Q_1$ of dimension 1 is the one-edge graph $K_2$, 
and a \textit{hypercube} $Q_n$ of dimension $n \geq 2$ is a Cartesian product $\Box_{i=1}^{n} K_2$. 
Isometric subgraphs of hypercubes are called \textit{partial cubes}. 
Median graphs form a subclass of partial cubes.
\textit{Djokovi\' c-Winkler relation} (briefly, \textit{relation $\Theta$}) plays an important role in the study of partial cubes.
Two edges $uv$ and $xy$ of a connected graph $G$ are said to be in \textit{relation $\Theta$},
denoted by $uv \Theta xy$, if $d_G(u, x) + d_G(v, y) \neq  d_G(u, y) + d_G(v, x)$.  
If a graph is isomorphic to an isometric subgraph of $Q_n$, we say that it can be {\em isometrically embedded} into $Q_n$.
The \textit{isometric dimension} of a partial cube $G$, denoted by $\mathrm{idim}(G)$,  
is the least integer $n$ for which $G$ embeds isometrically into a hypercube $Q_n$.
For any partial cube $G$, relation $\Theta$ is an equivalence relation on the edge set of $G$,
and $\mathrm{idim}(G)$ is the number of $\Theta$-classes of $G$ \cite{D73}.
{It is trivial that a partial cube of isometric dimension 0 is the one-vertex graph $K_1$.} 

Hypercubes {of dimension $n \ge 1$} can be defined equivalently using binary codes.
Let $\mathcal{B}^n=\{0, 1\}^n$ be the set of all \textit{binary codes} (or, \textit{binary strings}) of length $n \ge 1$.
Then a {\em hypercube} $Q_n$  has the vertex set $\mathcal{B}^n$,  
and two vertices of $Q_n$ are adjacent if the corresponding binary codes differ in precisely one position. 
A partial order $\leq$ can be defined on $\mathcal{B}^n$ with 
$u_1u_2 \ldots u_n \leq v_1v_2 \ldots v_n$ if $u_i \leq v_i$ for all $1 \le i \le n$,
{and the resulted poset is denoted by $(\mathcal{B}^n, \le)$.}
Daisy cubes were introduced as isometric subgraphs of hypercubes in \cite{KM19}. 
 A \textit{daisy cube}  is an induced subgraph of $Q_n$ generated by 
a nonempty poset {$(X, \le)$ contained in $(\mathcal{B}^n, \le)$}, which we denote by $Q_n(X)$ and define 
as $Q_n(X)=\langle \{ u \in \mathcal{B}^{n} \ | \ u \leq x \textrm{ for some } x \in X \} \rangle.$

The authors of the seminal paper \cite{KM19} observed the following results of daisy cubes. 
\begin{lemma} \cite{KM19} \label{L:DaisyCubeProperties}
Daisy cubes are partial cubes. Moreover, for a daisy cube $Q_n(X)$ with $X \subseteq \mathcal{B}^n$,
(i) if  $\widehat{X}$  is the set of maximal elements of the poset $(X, \leq)$, then  
$Q_n(X)=Q_n(\widehat{X}) = \left\langle \cup_{x \in \widehat{X}} I_{Q_n}(x, 0^n)  \right\rangle$, 
where $0^n$ denotes a string of $n$ zeros;
(ii) $\mathrm{idim}(Q_n(X))= \deg_{Q_n(X)}(0^n)$.
\end{lemma}

{Note that if $\widehat{X}$  is the set of maximal elements of the poset $(X, \leq)$ 
where $X \subseteq \mathcal{B}^n$, then 
$\widehat{X}$ forms an anti-chain of $(X, \leq)$, and thus of $(\mathcal{B}^n, \leq)$.
Hence,} if $G$ is a daisy cube, 
then there is a nonempty subset $\widehat{X} \subseteq \mathcal{B}^n$ for some positive integer $n$  
such that $\widehat{X}$ forms an anti-chain of $(\mathcal{B}^n, \leq)$ 
and $G=\langle \{ u \in \mathcal{B}^{n} \ | \ u \leq x \textrm{ for some } x \in \widehat{X} \} \rangle 
= \left\langle \cup_{x \in \widehat{X}} I_{Q_n}(x, 0^n)  \right\rangle$;
and we call such a nonempty subset $\widehat{X} \subseteq \mathcal{B}^n$ as the set
of \textit{maximal vertices} of $G$.

\subsection{Perfect matchings of plane bipartite graphs}\label{S2-PM}

A \textit{perfect matching} $M$ of a graph $G$ is a subset of  $E(G)$
such that every vertex of  $G$ is incident with exactly one edge from $M$. 
If $H$ is a path or a cycle of $G$, then $H$ is \textit{$M$-alternating} if  edges of $H$ are alternately  in $M$ and  out of $M$. 
A path $P$ of  $G$ is called a \textit{handle} if all internal vertices (if exist) of $P$ are degree-2 vertices of $G$,
and each end vertex of $P$ has degree at least three in $G$ \cite{CC13}. 
A handle with more than one edge is called a \textit{nontrivial handle}. 
It is clear that any nontrivial handle of $G$ is $M$-alternating for any perfect matching $M$ of $G$.

Assume that $G$ is a plane graph. Regions bounded by edges of $G$ are called faces of $G$.
 A face $s$ of $G$  is called a \textit{finite face} 
(or, an \textit{inner face}) if  $s$ represents  a finite region,  and \textit{the infinite face} otherwise.
In addition, we denote the set of edges enclosing a face $s$ of $G$ by $E(s)$. 
The subgraph induced by the edges in $E(s)$ is the \textit{periphery of $s$}. 
The periphery of the infinite face is also called the \textit{periphery of $G$}.
Vertices on the periphery of $G$ are called \textit{exterior vertices} and the remaining vertices are \textit{interior vertices}. 
Also, edges on the periphery of $G$ are called \textit{exterior edges} and the remaining edges are \textit{interior edges}. 
A handle of $G$ is called  an \textit{exterior handle}  (respectively, an \textit{interior handle}) if all of its edges are  exterior edges 
(respectively, interior edges) of the plane graph.  
An \textit{outerplane graph} is a plane graph  whose vertices are all exterior. 

Let $G$ be a plane bipartite graph with a perfect matching $M$. 
If the periphery of a face $f$ of $G$ is an $M$-alternating cycle, 
then $f$  is called {\textit{$M$-resonant} \cite{ZZ00},
and we say that $f$ is \textit{resonant} briefly  if there is no need to specify the perfect matching $M$.}
The Fries number was initially introduced for benzenoid  hydrocarbons \cite{F27}, 
and naturally extended to plane bipartite graphs, for example, see \cite{AA07}.
Let $G$ be a plane bipartite graph and $S$ be the set of finite faces of $G$ with $|S| = n$ for some positive integer $n$.
The  \textit{Fries number} of $G$ is the maximum cardinality of a subset $S' \subseteq S$ 
satisfying the property that there exists a perfect matching $M$ of $G$ such that each finite face in $S'$ is $M$-resonant.
If the Fries number of $G$ is $n$, then there exists a perfect matching $M$ of $G$ such that each finite face of $G$ is $M$-resonant.

An edge of a graph with a perfect matching 
is called \textit{allowed} if it is contained in some perfect matching of the graph, and \textit{forbidden} otherwise. 
Note that if a face $s$ is $M$-resonant for some perfect matching $M$ of a plane bipartite graph $G$, 
then each edge in $E(s)$  is allowed since  each edge in $E(s)$ is contained in either $M$ or $M \oplus E(s)$.
A graph is said to be \textit{elementary} if all its allowed edges form a connected subgraph. 
It was shown that a bipartite graph is elementary if and only if it is connected and each edge is allowed,
{and an elementary bipartite graph with more than two vertices is 2-connected} \cite{LP86}.
In addition, a connected plane bipartite graph with more than two vertices  
is elementary if and only if  each face (including the infinite face) is resonant \cite{ZZ00}.

For a plane bipartite graph $G$ with a perfect matching, the \textit{elementary components} of $G$ are  
the components of the  subgraph obtained from $G$ by removing all forbidden edges of $G$. 
Obviously,  elementary components of $G$ are plane elementary bipartite graphs. 
The concept of a plane weakly elementary bipartite graph was first introduced in \cite{ZZ00} for connected graphs.
It is more practical to extend the concept to include non-connected graphs too. 
A  plane bipartite graph (not necessarily connected) is called \textit{weakly elementary} 
if deleting all forbidden edges does not produce any new finite face. 
By definition, any plane elementary bipartite graph is also weakly elementary.
If a plane weakly elementary bipartite $G$ has a finite face $s$ which is not a finite face of any elementary component of $G$,
then $s$ contains a forbidden edge of $G$, and so $s$ cannot be $M$-resonant for any perfect matching $M$ of $G$.

\subsection{Resonance graphs}\label{S2-RG}

The \textit{resonance graph} (also called \textit{$Z$-transformation graph}) $R(G)$ of a plane bipartite graph $G$ is a graph 
whose vertices are the  perfect matchings of $G$, and two perfect matchings $M_1,M_2$ are adjacent 
whenever their symmetric difference $M_1 \oplus M_2$ forms exactly
one cycle that is the  periphery of some finite face $s$ of $G$ \cite{ZZY04}.
In this case, we say that the edge $M_1M_2$  {of $R(G)$} has the \textit{face-label} $s$.
It is clear that if $s$ is a face-label of an edge $M_1M_2$ in $R(G)$, 
then $s$ is both $M_1$-resonant and $M_2$-resonant, and each edge on the periphery of $s$ is allowed.
It is well known \cite{ZLS08}  that if $G$ is a plane weakly elementary bipartite graph, then 
its resonance graph $R(G)$ is a median graph.  Since any median graph is a partial cube,
the relation $\Theta$  is an equivalence relation on the edge set of $R(G)$.

Let $G$ be a plane weakly elementary bipartite graph with elementary components $G_1$, $G_2$, $\ldots, G_t$.
Then the resonance graph $R(G)$ is a Cartesian product of resonance graphs $R(G_i)$ for all $1 \le i \le t$, 
that is, $R(G)=\Box_{i=1}^{t}R(G_i)$, which is well known and can be explained as follows.
Note that any perfect matching $a$ of $G$ is a disjoint union of $a_1, a_2, \ldots, a_t$ 
where $a_i$ is the restriction of $a$ on $G_i$ for  all $1 \le i \le t$. 
Then each vertex of $R(G)$ can be represented by $a=(a_1, a_2, \ldots, a_t)$.
Two vertices $a=(a_1, a_2, \ldots, a_t)$ and $b=(b_1, b_2, \ldots, b_t)$ 
are adjacent in $R(G)$ if and only if the symmetric difference of $a$ and $b$ 
forms the periphery of a finite face $s$ of $G$. Since $G$ is weakly elementary, 
$s$ is a finite face of $G_i$ for some $1 \le i \le t$. Consequently, $a$ and $b$ are adjacent in $R(G)$ 
if and only if  the symmetric difference of $a$ and $b$ 
forms the periphery of a finite face $s$ of $G_i$ for some $1 \le i \le t$, which is equivalent to the fact that
 $a_i$ and $b_i$ are adjacent in $R(G_i)$ for some $1 \le i \le t$,
and $a_j=b_j$ for all $j \in \{1, 2, \ldots, t\} \setminus \{i\}$.
By the definition of a Cartesian product graph,  $R(G)$ is the Cartesian product $\Box_{i=1}^{t}R(G_i)$.
Note that if $G_i$ is $K_2$  for some $1 \le i \le t$, 
then $R(G_i)$ is the one-vertex graph, and contributes a trivial factor of  $\Box_{i=1}^{t}R(G_i)$.

\begin{lemma}\label{L:idim(R(G))=0}
Let $G$ be a plane bipartite graph with a perfect matching.
Then {the resonance graph} $R(G)$ is a daisy cube with $\mathrm{idim}(R(G)) = 0$ 
if and only if $G$ is a plane weakly elementary  bipartite graph
 such that each elementary component of $G$ is $K_2$.
\end{lemma}
\proof 
If {the resonance graph} $R(G)$ is a daisy cube with isometric dimension 0, {then} $R(G)$ is the one-vertex graph and therefore, 
graph $G$ has exactly one perfect matching. 
Hence, by removing all forbidden edges of $G$, 
the connected components of the obtained graph are all $K_2$. 
Consequently, $G$ is plane weakly elementary  bipartite graph such that each elementary component of $G$ is $K_2$.

Next, suppose that $G$ is a plane weakly elementary  bipartite graph
 such that each elementary component of $G$ is $K_2$. Note that the resonance graph of the one-edge graph $K_2$ is the one-vertex graph.
By the fact that the resonance graph $R(G)$ of a plane weakly elementary bipartite graph $G$ 
with the elementary components $G_i$ (where $1 \le i \le t$) is a Cartesian product $R(G)=\Box_{i=1}^{t}R(G_i)$,
we have that  $R(G)$ is the one-vertex graph, so $R(G)$ is a daisy cube with $\mathrm{idim}(R(G)) = 0$.
\qed\\

In this paper, we will focus on studying plane bipartite graphs whose resonance graphs are daisy cubes with 
isometric dimension at least 1. Note that any daisy cube is connected. We will need the following result in  \cite{ZZY04} 
which characterizes when the resonance graph of a plane bipartite graph is connected.

\vskip 0.1in

\begin{theorem} \cite{ZZY04} \label{con_res}
Let $G$ be a plane bipartite graph with a perfect matching. 
Then the resonance graph $R(G)$ is connected if and only if $G$ is weakly elementary.
\end{theorem}

By Lemma \ref{L:idim(R(G))=0} and Theorem \ref{con_res},
plane bipartite graphs whose resonance graphs are daisy cubes with 
isometric dimension at least 1 are weakly elementary and with at least one elementary component other than $K_2$.
We have seen that
each elementary component of a plane weakly elementary bipartite graph is  a plane elementary bipartite graph,
and any plane elementary bipartite graph with more than two vertices (that is, other than $K_2$) is 2-connected.

 The resonance graph of a catacondensed even ring system $G$ 
(respectively, a 2-connected outerplane bipartite graph $G$)
being a daisy cube was characterized in  \cite{br-tr-zi-1} (respectively,  \cite{br-tr-zi-2}) in terms of $G$ being ``regular''.
Since the terminology ``regular'' in \cite{br-tr-zi-1} and \cite{br-tr-zi-2} is different from that is 
usually adopted to define those graphs whose vertices all have the same degree,
we rename the concept using  ``angularly connected faces'' and ``linearly connected faces'' 
as a generalized version of the corresponding terminologies called
angularly connected hexagons and linearly connected hexagons, which were used in  \cite{Z18} to tell 
that the resonance graph of a kinky catacondensed benzenoid graph is a daisy cube.
Let $s$, $s'$, $s''$ be three finite faces of a 2-connected outerplane bipartite graph $G$ 
such that $s, s'$ have the common edge $e$, and
$s', s''$ have the common edge $f$. 
Then the triple $(s, s', s'')$ is called an \textit{adjacent triple of finite faces}.
Let $d_{L(G)}(e, f)$ be the distance between two edges $e$ and $f$ in the line graph $L(G)$ of $G$.
The adjacent triple of finite faces $(s, s', s'')$  is \textit{angularly connected} if $d_{L(G)}(e,f)$ is even, 
and \textit{linearly connected} otherwise.
One main result in \cite{br-tr-zi-1} can be rephrased as follows: {If} $G$ is a catacondensed even ring system, 
then its resonance graph $R(G)$ is a daisy cube  if and only if $G$ does not have any adjacent triple of finite faces
that is linearly connected.
The extension of the above result in \cite{br-tr-zi-2} for 
resonance graphs of $2$-connected outerplane bipartite graphs can be rephrased as follows.

\smallskip
\smallskip
\begin{theorem}\label{T:daisycube} \cite{{br-tr-zi-2}}
Let $G$ be a $2$-connected outerplane bipartite graph.
Then the resonance graph $R(G)$ is a daisy cube if and only if  $G$ does not have any adjacent triple of finite faces
that is linearly connected.
\end{theorem}

\section{Characterizations for plane elementary bipartite\\ graphs} 
\label{sec3}

In this section, we characterize plane elementary bipartite graphs whose resonance graphs are daisy cubes. 
It is known that \cite{LP86} that any plane elementary bipartite graph {other than $K_2$} is 2-connected.
We start with the following lemma which will be needed for the proof of our main theorem.

\begin{lemma} \label{L:InteriorHandle}
Let $G$ be a 2-connected plane bipartite graph with $n$ finite faces for some positive integer $n$. 
If the Fries number of $G$ is $n$,
then any interior vertex of $G$ has degree 2,  any exterior vertex of $G$ has degree at most 3,
and any  handle of $G$ has odd length. 
\end{lemma}

\proof If the Fries number of $G$ is $n$,
then $G$ has a perfect matching $M_A$ such that all finite faces of $G$ are $M_A$-resonant.

Suppose that $G$ has an interior vertex $u$ with $\deg_G(u)=t \ge 3$. 
Let $uv_1, uv_2, \ldots, uv_t$ be the edges of $G$ incident to $u$. 
Then all edges $uv_i$ ($1 \le i \le t$) are interior edges of $G$ since $u$ is an interior vertex of $G$,
and there are $t$ finite faces containing vertex $u$.
We can assume that $s_1, s_2, \ldots, s_t$ are the finite faces of $G$ containing vertex $u$
such that  $uv_i$ is a common edge of $s_{i-1}$ and $s_{i}$ for all $2 \le i \le t$, and $uv_1$ is a common edge of $s_t$ and $s_1$.
Note that exactly one of $uv_1, uv_2, \ldots, uv_t$ is contained in $M_A$. 
Without loss of generality, we can assume that $uv_1 \in M_A$.  
Then  $uv_2$ and $uv_3$  incident to $u$ are two consecutive edges on the periphery of $s_2$ and  not contained in $M_A$. 
So, $s_2$ is not $M_A$-resonant. This is a contradiction to the assumption. Therefore, any interior vertex of $G$ has degree 2.

Next, we show that  any exterior vertex of $G$ has degree at most 3.
Suppose that $G$ has an exterior vertex $w$ with $\deg_G(w)=d\ge 4$.  
Let $wz_1, wz_2, \ldots, wz_d$ be the edges of $G$ incident to $w$. 
Then all but two of these edges are interior edges of $G$ since $w$ is an exterior vertex of $G$,
and there are $d-1$ finite faces containing vertex $w$.
Without loss of generality, we can assume that $wz_1$ and $wz_d$ are exterior edges. 
Let $s_1, s_2, \ldots, s_{d-1}$ be the finite faces of $G$ containing vertex $w$
such that $wz_i$ is a common edge of $s_{i-1}$ and $s_{i}$ for all $2 \le i \le d-1$.
Note that exactly one of $wz_i$ where $1 \le i \le d$  is contained in $M_A$. 
There are two possibilities based on the edge $wz_i$ contained in $M_A$ is an exterior edge or not.
If the edge $wz_i$ contained in $M_A$ is an exterior edge, then we can assume that $wz_1 \in M_A$. 
It follows that  $wz_2$ and $wz_3$  incident to $w$ are two consecutive edges on the periphery of $s_2$ and not contained in $M_A$.
So, $s_2$ is not $M_A$-resonant.  This is a contradiction to the assumption.
If the edge $wz_i$ contained in $M_A$ is an interior edge, then we can assume that $wz_2 \in M_A$. 
It follows that  $wz_3$ and $wz_4$  incident to $w$ are two consecutive edges on the periphery of $s_3$ and not contained in $M_A$.
So, $s_3$ is not $M_A$-resonant. This is a contradiction to the assumption. Therefore, any exterior vertex of $G$ has degree at most 3.

We now show that each interior handle of $G$ has odd length. Suppose that $G$ has an interior handle $P$ with even length.
Then $P$ cannot be an edge and we can write $P=xx_1x_2\ldots x_{2t-1}y$ for some positive integer $t$.
Since all interior vertices $x_1, \ldots, x_{2t-1}$ of $P$ are degree-2 vertices of $G$,  it follows that
exactly one end edge of $P$ is not contained in $M_A$. 
Without loss of generality, we assume that $x_{2t-1}y$ is not contained in $M_A$.
Note that $y$ is an exterior vertex of $G$ with degree 3, 
so the two edges incident to $y$ (one of them is $x_{2t-1}y$)  
not contained in $M_A$  form two consecutive edges of a finite face $s$ of $G$.
Then $s$ is not $M_A$-resonant. This is a contradiction. 
Therefore, each interior handle of $G$ has odd length.

We further show that an interior handle of $G$ is either a single edge contained in $M_A$ 
or an $M_A$-alternating path with both end edges inside $M_A$.
If an interior handle of $G$ is an edge, then it is contained in $M_A$. 
Otherwise, similarly as above we can show that  $G$ has a finite face $s$ not $M_A$-resonant, which is a contradiction.
Assume that an interior handle of $G$ is nontrivial. Then it is
$M_A$-alternating since any internal vertex of the handle is a degree-2 vertex of $G$.
Moreover, two end edges of such an interior handle are either both inside $M_A$ or both outside $M_A$ since it has odd length.
Suppose that there is an interior handle whose two end edges are not contained in $M_A$. 
Then similarly we can show that $G$ has a finite face $s$ not $M_A$-resonant, which is a contradiction.
Therefore, each interior nontrivial handle of $G$  is $M_A$-alternating with both end edges inside $M_A$.

Finally, we consider exterior handles of $G$.
Any exterior nontrivial handle must be $M_A$-alternating.  
Obviously, any end vertex of an exterior handle is also an end vertex of some interior handle.
From the previous discussion on interior handles, it follows that for each exterior handle $P$ of $G$, 
either $P$ is a single edge not contained in $M_A$ or 
$P$ is $M_A$-alternating with no end edges contained in $M_A$. Therefore, any exterior handle of $G$ has odd length.
 \qed \\

We now introduce the following new terminology to describe the structure of plane elementary bipartite graphs 
whose resonance graphs are daisy cubes.

\begin{definition}\label{D:PCA}
Let $G$ be a plane elementary bipartite graph other than $K_2$. 
Then $G$ is called \textbf{peripherally 2-colorable} if  every vertex of $G$ has degree 2 or 3, 
vertices with degree 3 (if exist) are all exterior vertices of $G$,
and $G$ can be properly 2-colored black and white so that two vertices with the same color are nonadjacent,
and vertices with degree 3 (if exist) are alternatively black and white along the clockwise orientation of the periphery of $G$. 
\end{definition}

\begin{figure}[h!] 
\begin{center}
\includegraphics[scale=0.8, trim=0cm 0.5cm 0cm 0cm]{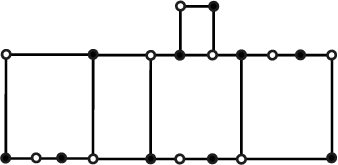}
\end{center}
\caption{\label{sl2} A peripherally 2-colorable  graph that is 2-connected outerplane bipartite.}
\end{figure}

In the next lemma, we demonstrate that if $G$ is a 2-connected outerplane bipartite graph,
then $G$ being peripherally 2-colorable coincides with $G$ being ``regular''  
from \cite{br-tr-zi-2}, which we rephrase as $G$ being ``without any adjacent triple of finite faces that is linearly connected''
in Section \ref{S2-RG}. See Figure \ref{sl2}.

\begin{lemma}\label{L:PCA-Regular}
Let $G$ be a 2-connected outerplane bipartite graph. {Then} $G$ is  peripherally 2-colorable
if and only if $G$  does not have any adjacent triple of finite faces that is linearly connected. 
\end{lemma} 
\proof  If $G$ has one finite face, then $G$ is an even cycle.
If $G$ has two finite faces, then they have  exactly one common edge since every vertex of $G$ 
is on the periphery of $G$. It follows that $G$ has two adjacent vertices of degree 3 and all other vertices
of $G$ have degree 2. Hence, the conclusion is trivial when $G$ has at most two finite faces.

Assume that $G$ is  {peripherally 2-colorable} and has at least three finite faces. 
Let $(s, s', s'')$ be an arbitrary adjacent triple of finite faces of $G$  
such that $e$ is the common edge of $s$ and $s'$, and $f$ is the common edge of $s'$ and $s''$.
Then the end vertices of $e$ and $f$ are degree-3 vertices of $G$ since every vertex of $G$ is at most 3.
Note that any exterior handle of $G$ has  odd number of edges
since vertices with degree 3 on the periphery of $G$ are alternately black and white 
 along the clockwise orientation of the periphery of $G$.
This implies that $d_{L(G)}(e, f)$, the distance between two edges $e$ and $f$ in the line graph of $G$, is even.  
Then $(s, s', s'')$ is  {angularly connected}.
Therefore, $G$  does not have any adjacent triple of finite faces that is linearly connected.

On the other hand, if $G$ does not have any adjacent triple of finite faces that is linearly connected, 
then any vertex of $G$ has degree at most 3.
Otherwise, if $G$ has a  vertex $w$ with degree at least 4,  
then there exists an adjacent triple of finite faces $(s, s', s'')$ containing $w$.
Assume that $e$ is the common edge of $s$ and $s'$, and $f$ is the common edge of $s'$ and $s''$.
Then $(s, s', s'')$ is linearly connected since $d_{L(G)}(e, f)=1$. This is a contradiction.
Hence,  any  vertex of $G$ has degree at most 3. 
{Recall that $G$ is 2-connected outerplane bipartite. 
Then} $G$ can be properly 2-colored  such that vertices with degree 3 are alternately black and white 
along the clockwise orientation of the periphery of $G$,
since any adjacent triple of finite faces $(s, s', s'')$  is angularly connected.
Therefore, $G$ is peripherally 2-colorable.
\qed\\

\begin{lemma}\label{L:InnerFace1-1ThetaClass}
Let $G$ be a plane elementary bipartite graph other than $K_2$.
If  the resonance graph $R(G)$ is a daisy cube, then
{$\mathrm{idim}(R(G))$ equals the number of finite faces of $G$,}
and  there is a 1--1 correspondence between the set of finite faces of  $G$ and the set of $\Theta$-classes of $R(G)$
{such that each $\Theta$-class of $R(G)$ has a unique face-label which is a finite face of $G$}.
\end{lemma}
 \proof Assume that $R(G)$ is a daisy cube with $\mathrm{idim}(R(G)) = n$ for some positive integer $n$.
It is trivial if $n=1$ since  $R(G)$ is the one-edge graph with one $\Theta$-class, 
and $G$ is an even cycle with one finite face.

Let $n >1$. Since $R(G)$ is a daisy cube with $\mathrm{idim}(R(G)) = n$, 
there is a subset $\widehat{X} \subseteq \mathcal{B}^n$ such that $\widehat{X}$ 
forms an anti-chain of $(\mathcal{B}^n, \leq)$ and 
$R(G)=\langle \{ u \in \mathcal{B}^{n} \ | \ u \leq x \textrm{ for some } x \in \widehat{X} \} \rangle 
= \left\langle \cup_{x \in \widehat{X}} I_{Q_n}(x, 0^n)  \right\rangle$; 
in particular notice that  $0^n$ is a vertex of $R(G)$. Furthermore, by the definition of a resonance graph, 
every edge $M_1M_2$ of $R(G)$ has a
face-label $s$ which is a finite face of $G$; 
in particular the periphery of $s$ is the symmetric difference of $M_1$ and $M_2$.
Hence, any two different edges incident to the same vertex  $0^n$ of $R(G)$ cannot have the same face-label.

By \cite{ZZ00}, a connected plane bipartite graph with more than two vertices is elementary if and only if 
each face is resonant.  Then any finite face $s$ of $G$  
is $M$-resonant for some perfect matching $M$ of $G$. 
It follows that $M'=M \oplus E(s)$ is another perfect matching of $G$, 
and $MM'$ is an edge of $R(G)$ with the face-label $s$. 
This implies that every finite face of  a plane elementary bipartite graph $G$ 
appears as a face-label of some edge in its resonance graph $R(G)$.
Recall that $R(G)= \langle \cup_{x \in \widehat{X}} I_{Q_n}(x, 0^n) \rangle$,
which is consisted of  hypercubes $\langle I_{Q_n}(x, 0^n) \rangle$ with a corner vertex $0^n$ for all $x \in \widehat{X}$.
It is well known \cite{ZOY09} that for a plane elementary bipartite graph $G$,  
antipodal edges of a 4-cycle in $R(G)$ have the same face-label.
Hence, for each  $x \in \widehat{X}$, face-labels of edges incident to $0^n$ in the hypercube $\langle I_{Q_n}(x, 0^n) \rangle$ 
are all possible face-labels appeared in the same hypercube.
It follows that face-labels of edges incident to $0^n$ in $R(G)$ form the set of all finite faces of $G$.

Now, we have shown that all face-labels of edges incident to $0^n$ are pairwise distinct 
and form  the set of finite faces of $G$. Hence, $\deg_{R(G)}(0^n)$ is the number of finite faces of $G$.
By Lemma \ref{L:DaisyCubeProperties} (ii), $\mathrm{idim}(R(G)) = \deg_{R(G)}(0^n)$.
{So, $\mathrm{idim}(R(G))$ equals the number of finite faces of $G$.}

{It is well known \cite{HIK11} that any two edges on a shortest path cannot be in the same $\Theta$-class
of a graph. Then any two edges incident to $0^n$ cannot be in the same $\Theta$-class of $R(G)$.
Recall that any daisy cube is a partial cube \cite{KM19}, and 
the isometric dimension of a partial cube is the number of $\Theta$-classes of the partial cube \cite{D73}.
By our assumption
that $R(G)$ is  daisy cube, we have that $R(G)$ is a partial cube, 
and $\mathrm{idim}(R(G))$ is the number of $\Theta$-classes of $R(G)$.
Recall that $\deg_{R(G)}(0^n)=\mathrm{idim}(R(G))$, 
and edges incident to $0^n$ are pairwisely contained in different $\Theta$-classes of $R(G)$.
Then each $\Theta$-class of $R(G)$ has exactly one edge incident to $0^n$.}

Hence, there is a 1--1 correspondence between the set of finite faces of $G$ and  the set of $\Theta$-classes of $R(G)$
{such that each $\Theta$-class of $R(G)$ has a unique face-label which is a finite face of $G$}.
\qed\\

\begin{figure}[h!] 
\begin{center}
\includegraphics[scale=0.8, trim=0cm 0.5cm 0cm 0cm]{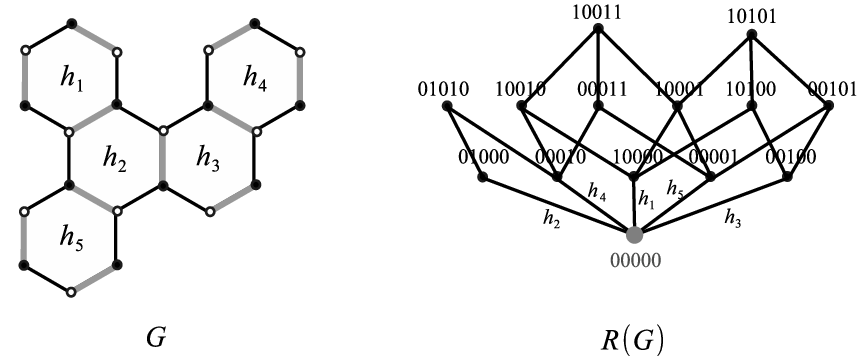}
\end{center}
\caption{\label{fig3} The resonance graph $R(G)$ of a kinky benzenoid graph $G$.}
\end{figure}

Lemma \ref{L:InnerFace1-1ThetaClass} can be visualized by the resonance graph of a kinky benzenoid graph in \cite{Z18}. 
By definition, we can see that a kinky benzenoid graph is a peripherally 2-colorable graph 
that is 2-connected outerplane bipartite.
{Edges of $G$ in Figure \ref{fig3} colored with grey form a perfect matching corresponding to the vertex of $R(G)$ with 
the binary code $00000$.}

The \textit{subdivision} of an edge $uv$ of a graph $G$ is an operation where we add a new degree-2 vertex $w$, 
and replace the edge $uv$ by two new edges $uw$ and $wv$ incident to $w$ while keeping previous neighbors 
of $u$ and $v$. The reverse operation, \textit{smoothing out} a degree-2 vertex $w$ of $G$ 
is to remove $w$, and replace two edges $uw$ and $wv$ incident to $w$  
with a new edge $uv$ while keeping neighbors of $u$ and $v$ different from $w$.

\begin{theorem}\label{T:R(G)DaisyCube}
Let $G$ be a plane elementary bipartite graph other than $K_2$. {Let $n$ be a positive integer.}
Then the resonance graph $R(G)$ is a daisy cube with $\mathrm{idim}(R(G)) = n$ 
if and only if the Fries number of $G$ is $n$, where $n$ is the number of finite faces of $G$.
\end{theorem}
\proof  
If $n=1$, then the conclusion is trivial since $G$ is an even cycle and $R(G)$ is an edge. Let $n>1$.

{\it Necessity:}  Assume that $R(G)$ is a daisy cube  with $\mathrm{idim}(R(G)) = n$.
Then by the definition of a daisy cube, we know that $0^n$ is a vertex of $R(G)$. 
 Let $M_A$ be the perfect matching of $G$ with the binary code representation $0^n$ as a vertex of $R(G)$. 
 By the proof of Lemma \ref{L:InnerFace1-1ThetaClass}, {$n$ is
 the number of finite faces of $G$, and} edges incident to $0^n$ have pairwise different face-labels
 which form  the set of $n$ finite faces of $G$. It follows that if $s$ is an arbitrary finite face of $G$,
then there is an edge $M_AM$  such that the symmetric difference of $M_A$ and $M$ is the periphery of $s$,
and so $s$ is $M_A$-resonant.  
Therefore,  every finite face of $G$ is $M_A$-resonant, and so the Fries number of $G$ is $n$.

{\it Sufficiency:}  Assume that the Fries number of $G$ is $n$, where $n$ is the number of finite faces of $G$.
Then $G$ has a  perfect matching $M_A$ such that all finite faces of $G$ are $M_A$-resonant. 
By  \cite{LP86},  a plane elementary bipartite graph {other than $K_2$} is 2-connected.
Therefore, a plane elementary bipartite graph $G$ with $n$ finite faces for some positive integer $n$ is 2-connected.
By Lemma \ref{L:InteriorHandle},  any interior vertex of $G$ has degree 2, 
any exterior vertex of $G$ has degree at most 3, and any handle of $G$ has odd length.
By definition, it follows that $G$ is peripherally 2-colorable. We distinguish two cases based whether $G$ 
has interior vertices or not.

Case 1. If $G$ has no interior vertices, then $G$ is a 2-connected outerplane bipartite graph
which is also peripherally 2-colorable.
By  Lemma \ref{L:PCA-Regular}, $G$ does not have any adjacent triple of finite faces that is linearly connected.
Hence, by Theorem \ref{T:daisycube},  $R(G)$ is a daisy cube. Moreover,  by Lemma \ref{L:InnerFace1-1ThetaClass}, 
$\mathrm{idim}(R(G))$ equals the number of finite faces of $G$ and therefore, $\mathrm{idim}(R(G)) =n$. The sufficiency is done.

Case 2. If $G$ has some interior vertices, then $G$ contains at least one interior nontrivial handle.
Since $G$ is peripherally 2-colorable, two end vertices of any interior handle have degree 3, so they are exterior vertices of $G$. 
Therefore, if an interior nontrivial handle has the property that its two end vertices $u$ and $v$ are adjacent in $G$, 
then the edge $uv$ must be an exterior edge of $G$.

Now for each interior nontrivial handle $P$ of $G$,  
we apply the following operations based on whether two end vertices $u,v$ of $P$ are adjacent in $G$ or not.

Subcase 2.1. If $u$ and $v$ are not adjacent in $G$, then replace $P$ by an interior edge $e'_P$ 
 by  smoothing out degree-2 vertices of $P$.

Subcase 2.2. If  $u$ and $v$ are adjacent in $G$, then we have observed above that $uv$ must be an exterior edge of $G$.
To avoid generating multiple edges, we replace the exterior edge $uv$ by an exterior nontrivial handle $P'_{uv}$ with odd length 
by  edge subdivisions. 
Then replace $P$ by an {interior} edge $e'_P$ by  smoothing out degree-2 vertices of $P$.

\begin{figure}[h!] 
\begin{center}
\includegraphics[scale=0.8, trim=0cm 0.5cm 0cm 0cm]{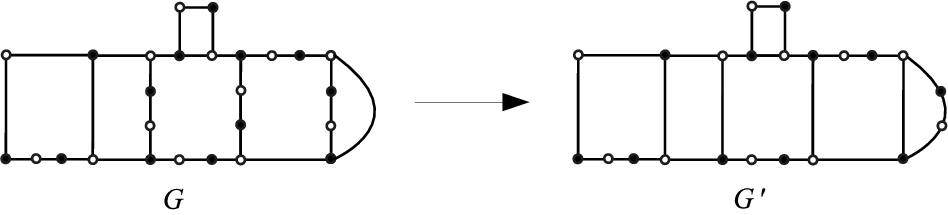}
\end{center}
\caption{\label{sl3} A peripherally 2-colorable graph $G$ is transformed into a 2-connected outerplane bipartite  graph $G'$
that is also peripherally 2-colorable by edge subdivisions and smoothings.}
\end{figure}

After above operations for each interior nontrivial handle $P$ of $G$, 
we obtain  {a 2-connected outerplane bipartite graph $G'$ 
which is also peripherally 2-colorable}. See Figure \ref{sl3}.

Similarly to the proof of Case 1,  we can show that $R(G')$ is a daisy cube. 
Moreover, the number of finite faces of $G'$ is the same as the number of finite faces of $G$ which is $n$. 
Again,  by Lemma \ref{L:InnerFace1-1ThetaClass}, 
$\mathrm{idim}(R(G'))$ equals the number of finite faces of $G'$ and therefore, $\mathrm{idim}(R(G')) =n$.

Note that there is an isomorphism between the inner dual of $G$ and the inner dual of $G'$
which maps a finite face $s_i$ of $G$ to a finite face $s'_i$ of $G'$ for all $1 \le i \le n$.
Moreover, there is a bijection $\phi$ between the set of perfect matchings of $G$ and the set of  perfect matchings of $G'$
such that $\phi$ maps a perfect matching $M$ of $G$ to a perfect matching $M'$ of $G'$ with the following properties.

(i)  If an interior nontrivial handle $P$ of $G$ is replaced by an interior edge $e'_P$ of $G'$, 
then $e'_P$ is contained in $M'$ of $G'$ if and only if two end edges of $P$ are contained in $M$ of $G$.

(ii) If an exterior edge $uv$ of $G$ is replaced by an exterior  nontrivial handle $P'_{uv}$  of $G'$ with odd length, 
then two end edges of $P'_{uv}$ are contained in $M'$ of $G'$ if and only if  $uv$ is contained in $M$ of $G$.

(iii) $M$ and $M'$ are identical on other unchanged handles and edges of $G$ and $G'$ during above operations.

It follows that for any two perfect matchings $M_1$ and $M_2$ of $G$,  
the symmetric difference $M_1 \oplus M_2$  is the periphery of a finite face $s_i$ of $G$
if and only if $\phi(M_1) \oplus \phi(M_2)$ is the periphery of a finite face $s'_i$ of $G'$ for some $1 \le i \le n$.
Therefore, the resonance graphs $R(G)$ and $R(G')$ are isomorphic. 
Consequently, $R(G)$ is also a daisy cube  with $\mathrm{idim}(R(G)) = n$.
\qed \\

We conclude this section a list of equivalent characterizations of plane elementary bipartite graphs 
whose resonance graphs are daisy cubes. 

\begin{corollary}\label{C:R(G)DaisyCube}
Let $G$ be a plane elementary bipartite graph other than $K_2$. 
{Let $n$ be a positive integer.}
Then the following statements are equivalent.

$(i)$ The resonance graph $R(G)$ is a daisy cube with $\mathrm{idim}(R(G))=n$.
 
$(ii)$  The Fries number of $G$ is $n$, where $n$ is the number of finite faces of $G$.
  
$(iii)$ $G$ is  peripherally 2-colorable and with $n$ finite faces.
\end{corollary}
\proof By Theorem \ref{T:R(G)DaisyCube}, $(i)$ and $(ii)$ are equivalent.
In the sufficiency part of the proof of Theorem \ref{T:R(G)DaisyCube} 
we showed that if the Fries number of $G$ is $n$, where $n$ is the number of finite faces of $G$, 
then  $G$ is  peripherally 2-colorable. This means $(ii)$ implies $(iii)$. 
In the same proof we also showed that  if $G$ is  peripherally 2-colorable and with $n$ finite faces, 
then the resonance graph $R(G)$ is a daisy cube with $\mathrm{idim}(R(G))=n$, 
so $(iii)$ implies $(i)$. Consequently, statements $(i)$, $(ii)$, and $(iii)$ are equivalent.
\qed\\

\section{Characterizations for plane bipartite graphs} 
\label{sec4}

In this section, we provide an answer for when the resonance graph of a plane bipartite graph is a daisy cube.
A daisy cube is a connected graph.
Thus by Theorem \ref{con_res},  if the resonance graph $R(G)$ of a plane bipartite graph $G$ is a daisy cube, 
then $G$ must be  weakly elementary. 
Let $G_1, G_2, \ldots, G_t$ be elementary components of a plane weakly elementary bipartite graph $G$ 
obtained by removing forbidden edges of $G$. 
Then $R(G)=R(G_1) \Box R(G_2) \Box \cdots \Box R(G_t)$ which is explained in detail in Section \ref{S2-RG}.

\begin{figure}[h!] 
\begin{center}
\includegraphics[scale=0.7, trim=0cm 0.5cm 0cm 0cm]{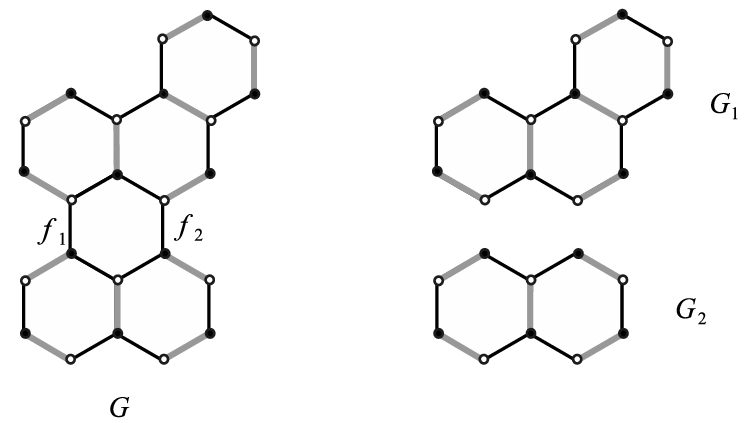}
\end{center}
\caption{\label{sl5} A plane bipartite graph $G$ {with} two elementary components $G_1$ and $G_2$.}
\end{figure}

\begin{figure}[h!] 
\begin{center}
\includegraphics[scale=0.7, trim=0cm 0.5cm 0cm 0cm]{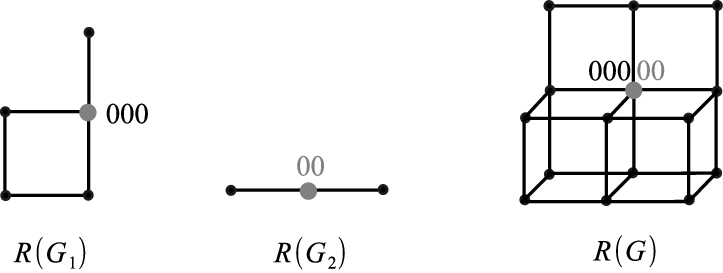}
\end{center}
\caption{\label{sl6} Resonance graphs of graphs $G_1$, $G_2$, and $G$ from Figure \ref{sl5}.}
\end{figure}

The plane bipartite graph $G$ in Figure \ref{sl5} has two elementary components $G_1$ and $G_2$ 
obtained after deleting two forbidden edges $f_1$ and $f_2$. 
It is easy to check that $G$ is weakly elementary, and $G_1, G_2$ are peripherally 2-colorable.
Moreover, $R(G)$ is the Cartesian product of $R(G_1)$ and $R(G_2)$, 
and all three resonance graphs are daisy cubes, see Figure \ref{sl6}.
{Edges of $G$, $G_1$, and $G_2$ colored with grey form perfect matchings corresponding to vertices of 
$R(G)$, $R(G_1)$, and $R(G_2)$ with the binary codes $00000$, $000$, and $00$, respectively.}

We will firstly show that a connected graph is a daisy cube if and only if each nontrivial Cartesian factor of the graph is a daisy cube.
The following terminologies will be used in the proof.
If $k \in \mathcal{B}^{n_1}$ and $l \in \mathcal{B}^{n_2}$ are binary codes,  
then the binary code obtained by concatenating $k$ and $l$ is denoted as $k \cdot l$ in $\mathcal{B}^{n_1+n_2}$. 
Moreover, if $k \in \mathcal{B}^n$ is a binary code and $n_1,n_2 \in \{ 1, \ldots, n \}$, 
then we use $k_{n_1}$ (respectively, $k_{\overline{n_2}}$) to denote the binary code  obtained from $k$ 
by taking the first $n_1$ positions of $k$ (respectively, the last $n_2$ positions of $k$). 
A \textit{$G$-layer} of $G \Box H$, denoted by $G \Box \{ h \}$, 
is a subgraph of $G \Box H$ induced by the vertices from the set $\{ (u,h) \ | \ u \in V(G) \}$,
where $h$ is a fixed vertex of $H$. 
 Similarly,  we can define a \textit{$H$-layer} $\{g \} \Box H$ for a fixed vertex $g$ of $G$.

\begin{theorem}\label{T:CartesianProduct-1}
{Let $n$ be a positive integer.} 
Then a Cartesian product $G \Box H$ of nontrivial graphs $G$ and $H$ 
is a daisy cube with $\mathrm{idim}(G \Box H)=n$  if and only if
$G$ and $H$ are daisy cubes with $\mathrm{idim}(G)=n_1$ and $\mathrm{idim}(H)=n_2$,
{where $n_1$ and $n_2$ are positive integers and} $n=n_1+n_2$.
\end{theorem}

\proof \textit{Sufficiency:}  This has been already observed on page 236 in \cite{EKM23}. 
For the sake of completeness, we provide a detail proof.
Suppose that $G$ and $H$ are daisy cubes with $\mathrm{idim}(G)=n_1$ and $\mathrm{idim}(H)=n_2$
{for some positive integers $n_1$ and $n_2$}.
Then let  ${X} \subseteq  \mathcal{B}^{n_1}$ be the set of maximal vertices of $G$ 
and ${Y}  \subseteq  \mathcal{B}^{n_2}$ be  the set of maximal vertices of $H$. 
{By Lemma \ref{L:DaisyCubeProperties} (i),} we can write
\begin{eqnarray*}
G&=& \langle \{ s \in \mathcal{B}^{n_1} \ | \ s \leq x \textrm{ for some } x \in {X} \} \rangle, \\
H&=& \langle \{ t \in \mathcal{B}^{n_2} \ | \ t \leq y \textrm{ for some } y \in {Y} \} \rangle.
\end{eqnarray*}
Let $U=\{s \cdot t \in \mathcal{B}^{n_1+n_2} \ | \ s \leq x \textrm{ for some } x \in X \textrm{ and } t \leq y \textrm{ for some } y \in Y \} $.
Then by slightly abusing notation we can write $U = V(G) \times V(H)$, 
where an ordered pair $(s,t)$ of codes  $s \in \mathcal{B}^{n_1}$, $t \in \mathcal{B}^{n_2}$ is identified with the code $s \cdot t$.
Note that $s_1 \cdot t_1$ and $s_2 \cdot t_2$ where $s_1,s_2 \in V(G)$ and $t_1,t_2 \in V(H)$ are two vertices of $G \Box H$.
Then $s_1 \cdot t_1$ and $s_2 \cdot t_2$ are adjacent in $G \Box H$ if and only if  either $s_1=s_2$ and $t_1$ and $t_2$  differ in exactly one position 
or $t_1=t_2$ and $s_1$ and $s_2$  differ in exactly one position
if and only if $s_1 \cdot t_1$ and $s_2 \cdot t_2$ differ in exactly one position.
It follows that $$G \Box H=\langle U \rangle =
\langle \{s \cdot t \in \mathcal{B}^{n_1+n_2} \ | \ s \leq x \textrm{ for some } x \in X, \textrm{ and } t \leq y \textrm{ for some } y \in Y \} \rangle.$$
Therefore,  $G \Box H$ is a daisy cube  {with $\mathrm{idim}(G \Box H)=n_1 + n_2$},
and {the set of maximal vertices} $M = \{x \cdot y \in \mathcal{B}^{n_1+n_2} \ | x \in X, y \in Y\}.$

\textit{Necessity:} Suppose that $G \Box H$ is a daisy cube with $\mathrm{idim}(G \Box H)=n$ {for some positive integer $n$}, 
where $G$ and $H$ are nontrivial graphs.
Then  $G \Box H$ is connected. So, $G$ and $H$ are connected graphs. 
Let $M$ be the set of maximal vertices of $G \Box H$.
Then $G \Box H= \langle \{ k \in \mathcal{B}^{n} \ | \ k \leq m \textrm{ for some } m \in M \} \rangle$.

Note that each edge of $G \Box H$ is contained either 
in a $G$-layer $G \Box \{h\}$ or in a $H$-layer $\{g\} \Box H$, where $g \in V(G)$ and $h \in V(H)$.
Next, we show that an edge from a $G$-layer and an edge from an $H$-layer of $G \Box H$ cannot be in relation $\Theta$. 
Assume that an edge $e=(u_1,v)(u_2,v)$ is from a $G$-layer $G \Box \{v\}$,
and an edge $f=(u,v_1)(u,v_2)$ is from an $H$-layer  $\{u\} \Box H$. 
By using Equation \eqref{en_dis},  we have
\begin{eqnarray*}
&&d_{G \Box H}((u_1,v),(u,v_1)) + d_{G \Box H}((u_2,v),(u,v_2)) \\
&=&  \left(d_{G}(u_1, u) + d_{H}(v, v_1) \right) +  \left( d_{G}(u_2,u) + d_{H}(v, v_2)\right)  \\
&=&  \left(d_{G}(u_1, u) + d_{H}(v, v_2)  \right) +  \left( d_{G}(u_2,u) + d_{H}(v, v_1) \right)  \\
&=&  d_{G \Box H}((u_1,v),(u,v_2)) + d_{G \Box H}((u_2,v),(u,v_1)).
\end{eqnarray*}
So, $e$ and $f$ are not in relation $\Theta$.
Therefore, the edges of a $\Theta$-class in $G \Box H$ are either all contained in $G$-layers or all contained in $H$-layers.

Note that $G \Box H$ has $n$ $\Theta$-classes as a partial cube of $Q_n$, 
and the vertex set $V(G \Box H) \subseteq \mathcal{B}^{n}$.
Then we can let $E_i$ be a $\Theta$-class of $G \Box H$ 
containing all edges between vertices differing in exactly position $i$ for some $1 \le i \le n$.
It follows that the $\Theta$-classes of $G \Box H$ are $E_1, \ldots, E_n$.
We have shown that the edges of a $\Theta$-class in $G \Box H$ are either all contained in $G$-layers or all contained in $H$-layers.  
Without loss of generality, we can assume that $\Theta$-classes  $E_1, \ldots, E_{n_1}$ are all contained in $G$-layers, 
and $\Theta$-classes $E_{n_1+1}, \ldots, E_{n}$ are all contained in $H$-layers. Let $n_2=n - n_1$.
{Note that both $n_1$ and $n_2$ are some positive integers since $G$ and $H$ are nontrivial graphs.} 
Then $V(G)$ and $V(H)$ can be represented using binary codes as follows:  
$V(G)=\{k_{n_1} \mid k \in V(G \Box H) \} \subseteq \mathcal{B}^{n_1}$ and 
$V(H)=\{k_{\overline{n_2}} \mid k \in V(G \Box H) \} \subseteq \mathcal{B}^{n_2}$.

Recall that $M$ is the set of maximal vertices of $G \Box H$. 
By the above conclusion, we can see that a vertex of $M$ 
can be written as $x \cdot y$ where $x \in  \mathcal{B}^{n_1}$ is a vertex of $G$, 
and $y \in  \mathcal{B}^{n_2}$ is a vertex of $H$.
Let $X=\{w_{n_1} \in V(G) \mid w \in M\}$ and $Y = \{w_{\overline{n_2}} \in V(H) \mid w \in M\}$.
It follows that 
\begin{eqnarray*}
G &=&\langle\{ s \in \mathcal{B}^{n_1} \ | \ s \leq x\textrm{ for some } x \in X \} \rangle \textrm{ and} \\
H &=& \langle \{ t \in \mathcal{B}^{n_2} \ | \ t \leq  y \textrm{ for some } y \in Y \}\rangle.
\end{eqnarray*}
Therefore, $G$ is a daisy cube  {with $\mathrm{idim}(G)=n_1$} and  the set of maximal vertices $X$, 
and $H$ is a daisy cube  {with $\mathrm{idim}(H)=n_2$} and the set of maximal vertices $Y$.
 \qed \\

\noindent
It is easy to show by induction that the previous theorem can be generalized to a daisy cube with more than two factors.
 
\begin{corollary}\label{C:CartesianProduct}
Let $R=R_1 \Box R_2 \Box \cdots \Box R_t$ be the Cartesian product of  nontrivial graphs $R_1, R_2, \ldots, R_t$.
{Let $n$ be a positive integer.} Then $R$ is a daisy cube with $\mathrm{idim}(R)=n$  
if and only if each $R_i$ is a daisy cube with $\mathrm{idim}(R_i)=n_i$, {where $n_i$ is 
a positive integer for all $1 \le i \le t$, and} $n=n_1+n_2 + \cdots + n_t$.
\end{corollary}

Finally, we conclude this section with {generalizations of our main results}.

\begin{theorem}\label{T:GeneralDaisyCube}
 Let G be a plane bipartite graph with a perfect matching. {Let $n$ be a positive integer.}
Then the resonance graph $R(G)$ is a daisy cube  with $\mathrm{idim}(R(G))=n$ 
if and only if $G$ is a plane weakly elementary bipartite graph, and for each of its elementary components $G_i$ other than $K_2$,
 the resonance graph $R(G_i)$ is a daisy cube with $\mathrm{idim}(R(G_i))=n_i$,  
 where {$n_i $ is a positive integer equal to} the number of finite faces of $G_i$ 
for all $1 \le i \le t$, and $n=n_1+n_2 + \cdots + n_t$.
\end{theorem}
\proof If $G$ is elementary, then  it is trivial by Theorem \ref{T:R(G)DaisyCube}. 
Assume that $G$ is not elementary. Then by Theorem \ref{con_res}, 
the resonance graph $R(G)$ is connected if and only if $G$ is plane weakly elementary bipartite.

If $R(G)$ is a daisy cube, then $R(G)$ is connected and so $G$ is plane weakly elementary bipartite. 
Let $G_1, G_2, \ldots, G_t$ be  {the} elementary components  {of $G$}
obtained by removing all forbidden edges of $G$. Then $t \ge 2$ by our assumption that $G$ is not elementary.
By the detail explanation in Section \ref{S2-RG}, we can see 
that $R(G)=R(G_1) \Box R(G_2) \Box \cdots \Box R(G_t)$.
Note that if $G_i$ is $K_2$ for some $1 \le i  \le t$, then $R(G_i)$ is the one-vertex graph, and contributes a trivial factor of $R(G)$.
Without loss of generality, we can assume that $G_i$ is different from $K_2$ for all $1 \le i \le t$, and so
$R(G_i)$ is a nontrivial factor of $R(G)$  for all $1 \le i \le t$.
Since $R(G)$ is a daisy cube, the conclusion follows  by 
{Theorem \ref{T:R(G)DaisyCube} and} Corollary \ref{C:CartesianProduct}.

On the other hand, if $G$
is a plane weakly elementary bipartite graph whose each elementary component $G_i$ different from $K_2$ has the property
that $R(G_i)$ is a daisy cube  with  $\mathrm{idim}(R(G_i))=n_i$,  where  {$n_i $ is a positive integer equal to} 
the number of finite faces of $G_i$
 for all $1 \le i \le t$, then by the fact that $R(G)=R(G_1) \Box R(G_2) \Box \cdots \Box R(G_t)$ and Corollary \ref{C:CartesianProduct}, 
 it follows that $R(G)$ is a daisy cube with $\mathrm{idim}(R(G))=n_1+n_2 + \cdots + n_t$, {which is a positive integer.}
\qed\\

By Theorem \ref{T:GeneralDaisyCube} and Corollary \ref{C:R(G)DaisyCube},  
we can obtain a list of equivalent characterizations of  
a plane bipartite graph whose resonance graph is a daisy cube,
which is a generalization of Corollary \ref{C:R(G)DaisyCube}. 
\begin{corollary}
Let $G$ be a plane bipartite graph with a perfect matching. 
{Let $n$ be a positive integer.}
Then the following statements are equivalent:

(i) The resonance graph $R(G)$ is a daisy cube with $\mathrm{idim}(R(G))=n$.

(ii) $G$ is plane weakly elementary bipartite, and for each of its elementary components $G_i$ other than $K_2$,
 the Fries number of $G_i$ is $n_i$,  
 where {$n_i$ is a positive integer equal to} the number of finite faces of $G_i$  for all $1 \le i \le t$, and $n=n_1+n_2 + \cdots + n_t$.
 
(iii)  $G$ is plane weakly elementary bipartite,  and for each of its elementary components $G_i$ other than $K_2$, 
$G_i$ is peripherally 2-colorable and with $n_i$  finite faces, where {$n_i$ is a positive integer for all  $1 \le i \le t$, 
and} $n=n_1+n_2 + \cdots + n_t$.
\end{corollary}

\smallskip

\noindent{\bf Acknowledgment:} Simon Brezovnik, Niko Tratnik, and Petra \v Zigert Pleter\v sek acknowledge the financial support from the Slovenian Research and Innovation Agency: research programme No.\ P1-0297 (Simon Brezovnik, Niko Tratnik, Petra \v Zigert Pleter\v sek), projects No.\ J1-4031, J2-2512 (Simon Brezovnik), N1-0285 (Niko Tratnik), and L7-4494 (Petra \v Zigert Pleter\v sek).  All
four authors thank the Slovenian Research and Innovation Agency for financing our bilateral project between Slovenia and the USA (title: \textit{Structural properties of resonance graphs and related concepts}, project No. BI-US/22-24-158).
  
 {The authors would like to thank the referees for their helpful comments.}


\begin{thebibliography}{99}
\bibitem{AA07} 
H. Abeledo, G. W. Atkinson, 
Unimodularity of the Clar number problem, 
Linear Algebra Appl. 420 (2007) 441--448.
 
\bibitem{br-tr-zi-1} 
S. Brezovnik, N. Tratnik, P. \v Zigert Pleter\v sek, 
Resonance graphs of catacondensed even ring systems, 
Appl. Math. Comput. 374 (2020), 125064, 9 pp.

\bibitem{br-tr-zi-2} 
S. Brezovnik, N. Tratnik, P. \v Zigert Pleter\v sek, 
Resonance graphs and a binary coding of perfect matchings of outerplane bipartite graphs, 
MATCH Commun. Math. Comput. Chem.  90 (2023) 453--468.

\bibitem{BCTZ23} 
S. Brezovnik, Z. Che, N. Tratnik, P. \v Zigert Pleter\v sek, 
Outerplane bipartite graphs with isomorphic resonance graphs, 
Discrete Appl. Math. 343 (2024) 340--349.

\bibitem{C18}
Z. Che,  
Structural properties of resonance graphs of plane elementary bipartite graphs, 
Discrete Appl. Math. 247 (2018) 102--110.

\bibitem{C19}
Z. Che, 
A characterization of the resonance graph of an outerplane bipartite graph,
Discrete Appl. Math. 258 (2019) 264--268.


\bibitem{C21} 
Z. Che, 
Peripheral convex expansions of resonance graphs,
Order 38 (2021) 365--376.  


\bibitem{CC13}
Z. Che, Z. Chen, 
Forcing faces in plane bipartite graphs (II), 
Discrete Appl. Math. 161 (2013) 71--80. 


\bibitem{D73}
D. \v{Z}. Djokovi\'c,
Distance preserving subgraphs of hypercubes,
J. Combin. Theory Ser. B  14 (1973) 263--267.

\bibitem{el-basil-93/1}
S.~El-Basil, 
Kekul\'e structures as graph generators, 
J. Math. Chem. 14 (1993) 305--318.

\bibitem{el-basil-93/2}
S.~El-Basil, 
Generation of lattice graphs. An equivalence relation  on Kekul\'e counts of catacondensed benzenoid hydrocarbons,
J. Mol. Struct.: THEOCHEM 288 (1993) 67--84.

\bibitem{EKM23}
\"{O}. E\v{g}ecio\v{g}lu, S. Klav\v{z}ar and M. Mollard,
\emph{Fibonacci cubes with applications and variations},
World Scientific Publishing Co Pte Ltd, 2023.

\bibitem{F27}
K. Fries, 
Uber Byclische Verbindungen und ihren Vergleich mit dem Naphtalin, 
Ann. Chem. 454 (1927) 121--324.

\bibitem{grundler-82}
W.~Gr{\"u}ndler, 
Signifikante Elektronenstrukturen fur benzenoide  Kohlenwasserstoffe, 
Wiss. Z. Univ. Halle 31 (1982) 97--116.

\bibitem{HIK11}
R. Hammack, W.  Imrich, S. Klav\v{z}ar,  
\emph{Handbook of product graphs.}  Second edition.  
Discrete Mathematics and its Applications (Boca Raton). CRC Press, Boca Raton, FL, 2011.


\bibitem{KM19}
S. Klav\v zar, M. Mollard, 
Daisy cubes and distance cube polynomial,
European J. Combin. 80 (2019) 214--223.


\bibitem{KZ05} 
S. Klav\v zar, P. \v Zigert Pleter\v sek, 
Fibonacci cubes are the resonance graphs of fibonaccenes, 
Fibonacci Quart.  43 (2005) 269--276.


\bibitem{LZ03}
P. C. B. Lam, H. Zhang, 
A distributive lattice on the set of perfect matchings of a plane bipartite graph, 
Order 20 (2003) 13--29. 

\bibitem{LP86}
L. Lovasz, M. D. Plummer, 
\emph{Matching theory}, North-Holland
Mathematics Studies, 121. Annals of Discrete Mathematics, Vol 29,
North-Holland Publishing Co., Amsterdam, 1986.

\bibitem{SKG06} 
K. Salem, S. Klav\v zar, I. Gutman, 
On the role of hypercubes in the resonance graphs of benzenoid graphs, 
Discrete Math. 306 (2006) 699--704.

\bibitem{SKVZ09} 
K. Salem, S. Klav\v zar, A. Vesel, P. \v Zigert, 
The Clar formulas of a benzenoid system and the resonance graph, 
Discr. Appl. Math. 157 (2009) 2565--2569.

\bibitem{T20}
A. Taranenko, 
Daisy cubes: A characterization and a generalization, 
European J. Combin. 85 (2020), 103058, 10 pp.
 

\bibitem{TV12}
A. Taranenko, A. Vesel, 
1-Factors and characterization of reducible faces of plane elementary bipartite graphs,
Discuss. Math. Graph Theory 32 (2012) 289--297.


\bibitem{TZ12}
A. Taranenko, P. \v Zigert Pleter\v sek,
Resonant sets of benzenoid graphs and hypercubes of their resonance graphs,
MATCH Commun. Math. Comput. Chem. 68 (2012) 65--77.

\bibitem{TD23} 
N. Tratnik, D. Ye, 
Resonance graphs on perfect matchings of graphs on surfaces, 
Graphs Combin. 39 (2023) Paper No. 68, 15 pp.

\bibitem{V21} 
A. Vesel, 
Efficient proper embedding of a daisy cube, 
Ars Math. Contemp. 21 (2021) Paper No. 7, 12 pp.

\bibitem{zhgu-88}
F. Zhang, X. Guo, R. Chen, 
$Z$-transformation graph of perfect matching of hexagonal systems, 
Discrete Math. 72 (1988) 405--415.

\bibitem{Z06}
H. Zhang, 
$Z$-transformation graphs of perfect matchings of plane bipartite graphs: a survey,
MATCH Commun. Math. Comput. Chem. 56 (2006) 457--476.


\bibitem{ZLS08}
H. Zhang, P. C. B. Lam, W. C. Shiu, 
Resonance graphs and a binary coding for the 1-factors of benzenoid systems, 
SIAM J. Discrete Math. 22 (2008) 971--984.


\bibitem{ZOY09}
H. Zhang, L. Ou and H. Yao,
Fibonacci-like cubes as $Z$-transformation graphs,
Discrete Math. 309 (2009)1284--1293. 

\bibitem{ZZ00}
H. Zhang, F. Zhang, 
Plane elementary bipartite graphs,
Discrete Appl. Math. 105  (2000) 291--311.

\bibitem{ZZY04}
H. Zhang, F. Zhang, H. Yao, 
$Z$-transformation graphs of perfect matchings of plane bipartite graphs, 
Discrete Math. 276  (2004) 393--404.

\bibitem{Z18} 
P. \v Zigert Pleter\v sek, 
Resonance graphs of kinky benzenoid systems are daisy cubes, 
MATCH Commun. Math. Comput. Chem. 80 (2018)  207--214.
\end{thebibliography}
\end{document}